\documentclass[11pt]{amsart}
\usepackage[english]{babel}
\usepackage{graphicx,amsmath,amssymb,color, enumerate}
\numberwithin{equation}{section}
\usepackage{enumerate}
\usepackage{color}

\newtheorem{theorem}{Theorem}[section]

\newtheorem{remark}[theorem]{Remark}

\newtheorem{proposition}[theorem]{Proposition}

\newcommand{\ee}{\mathrm{e}}
\newcommand{\dd}{\mathrm{d}}

\newcommand{\Dom}{\mathsf{Dom}}

\newcommand{\N}{\mathbb{N}}
\newcommand{\Z}{\mathbb{Z}}
\newcommand{\R}{\mathbb{R}}

\newcommand{\Ab}{\mathbf{A}}
\newcommand{\sA}{\mathsf A}
\newcommand{\sS}{\mathsf S}

\newcommand{\Dir}{\mathrm{Dir}}
\newcommand{\Neu}{\mathrm{Neu}}
\newcommand{\sH}{\mathsf{H}}
\newcommand{\fh}{\mathfrak h}

\newcommand{\sJ}{\mathsf{J}}
\newcommand{\spec}{\mathsf{sp}}

\newcommand{\ii}{\,\mathrm{i}}

\DeclareMathOperator{\curl}{curl}

\begin{document}

\title[Counting eigenvalues below the lowest Landau level]{Counting eigenvalues below the lowest Landau level}

\author{S{\o}ren Fournais}
\address[S. Fournais]{Department of Mathematical Sciences,
University of Copenhagen,
Universitetsparken 5,
DK-2100 Copenhagen OE,
Denmark}
\email{fournais@math.ku.dk}

\author{Ayman Kachmar}
\address[A. Kachmar]{The Chinese University of Hong Kong, Shenzhen, Guangdong, 518172, P.R. 
China.}
\email{akachmar@cuhk.edu.cn}

\date{\today}

\begin{abstract}
For the magnetic Laplacian on a bounded planar domain, imposing Neumann boundary conditions produces eigenvalues below the lowest Landau level. If the domain has two boundary components and one imposes a Neumann condition on one component and a Dirichlet condition on the other, one gets fewer such eigenvalues than when imposing Neumann boundary conditions on the two components. We quantify this observation for two models: the  strip and the annulus. In both models one can separate variables and deal with a family of fiber operators, thereby reducing the problem to counting band functions, the eigenvalues of the fiber operators. 
\end{abstract}

\maketitle

\section{Introduction}

\subsection{The magnetic Laplacian in multiply connected domains}
Let $\Omega\subset \R^2$ be open, bounded and with smooth boundary $\Gamma$ consisting of $N+1$ simple closed curves $\Gamma_0,\cdots,\Gamma_N$, where we assume that $\Gamma_1,\cdots,\Gamma_N$ are contained in the simply connected set $\Omega_0$ with boundary $\Gamma_0$. Consider the magnetic Laplacian 
\[\sH_\Omega=(-\ii h\nabla-\Ab)^2\]
in $L^2(\Omega)$, where $h>0$ is a small parameter (the semi-classical parameter) and the vector potential $\Ab:\Omega\to\R^2$ is  smooth and satisfies (constant magnetic field)
\[\curl\Ab=1.\]

We can obtain various self-adjoint realizations of the operator $\sH_\Omega$ by imposing Dirichlet or Neumann boundary conditions on the different boundary components. If we impose the Dirichlet boundary condition, $u=0$, on all of $\partial \Omega$, we get that the infimum of the spectrum is larger than the lowest Landau level $h$. However, if we impose a Neumann boundary condition on at least one boundary component, which is the case we will deal with in this paper, then eigenvalues below $h$ will appear (at least for $h$ small enough, see below). Our aim in this paper is to count these eigenvalues in the `semi-classical' limit, $h\to0$.

To count these eigenvalues in general has surprisingly enough turned out to be difficult. This is in part due to the fact that the semi-classical Weyl term is discontinuous at the Landau levels $(2n+1)h$, $n\in {\mathbb N} \cup \{0\}$. 
Also, for most of the eigenvalues, the effect of the boundary is exponentially small, making it delicate to determine the sign of such a small effect.
However, for the specific case of the lowest Landau level $h$, the recent paper \cite{FFGKS} succeeds in proving the correct formula in the specific case of purely Neumann boundary conditions or more general Robin conditions. 
The number of eigenvalues below $h$ in this case is proven to be
$$ \sim (2\pi h)^{-1}|\Omega|,$$
to leading order as $h \rightarrow 0$.

A Dirichlet condition formally corresponds to taking the Robin parameter to infinity. Therefore a simple intuition, based on the strong localisation of magnetic eigenfunctions, would suggest that 
in the formula for the counting function
in the case of mixed Dirichlet-Neumann boundary conditions, the area $|\Omega|$ should be replaced by the area of the smaller set where the nearest boundary point has a Neumann (instead of a Dirichlet) condition.
We will below give two concrete geometric examples of mixed magnetic Dirichlet-Neumann boundary conditions. In the first, the prediction of the simple intuition described above is verified. However, in the second case - the case of an annulus - the asymptotic result disagrees with the intuition and therefore points to a more involved geometric picture for the counting function.

\subsection{Mixed Dirichlet-Neumann boundary conditions}
Let us introduce the Neumann boundary condition for the operator $\sH_\Omega$. We denote by $\nu_0,\cdots,\nu_N$ the unit normal vectors of $\Gamma_0,\cdots,\Gamma_N$ respectively, and we assume that they point inward into $\Omega$. By a Neumann boundary condition on the $k$'th component $\Gamma_k$, we mean
\begin{equation}\label{eq:N-bc-i}
\nu_k\cdot (-\ii h\nabla -\Ab)u=0\quad\mbox{ on }\Gamma_k.
\end{equation}
Consider  non-empty, disjoint sets $K^{\Dir}$ and $K^{\Neu}$ such that $K^{\Dir}\cup K^{\Neu}=K_N:=\{0,1,\cdots,N\}$, and let us assume that $\sH_\Omega$ is the self-adjoint operator in $L^2(\Omega)$ obtained by imposing
\begin{enumerate}[i)] 
\item
Dirichlet boundary conditions, $u=0$, on every boundary component $\Gamma_k$ with $k\in K^{\Dir}$;
\item Neumann boundary condition, as in \eqref{eq:N-bc-i}, on every boundary component $\Gamma_k$ with
$k\in K^{\Neu}$.
\end{enumerate}
These boundary conditions serve to realize $\sH_\Omega$ as a positive self-adjoint operator in $L^2(\Omega)$ with compact resolvent and therefore spectrum consisting of eigenvalues $\{ \lambda_j(h,\Omega,\Ab) \}_{j=1}^{\infty}$ (numbered in increasing order including multiplicity) with $\lambda_j(h,\Omega,\Ab) \rightarrow \infty$ as $j \rightarrow \infty$.

More generally,  the above definition of the operator still makes sense if $\Omega$ is a  2-manifold with boundary embedded in $\R^3$, for instance the strip (or cylinder) $\mathbb S^1\times(0,L)$. We will discuss this particular example more in detail in Section~\ref{sec:strip}  later on.

\subsection{Earlier results}
Given $\lambda\in(0,1]$ and $h>0$,  we introduce the set
\begin{equation}\label{eq:def-set}
\sJ(\lambda h,\Omega,\Ab):=\{j\colon \lambda_j(h,\Omega,\Ab)<\lambda h \},
\end{equation}
and we denote by
\[N(\lambda h,\Omega,\Ab)=|\sJ(\lambda h,\Omega,\Ab)| \]
the cardinal of $\sJ(\lambda h,\Omega,\Ab)$, i.e. the number of eigenvalues of $\sH_\Omega$ below $\lambda h$, counting multiplicity.

Quite recently, it was proved in \cite{FFGKS} that for simply connected domains and Neumann boundary condition
\begin{equation}\label{eq:nb-sc}
N(h,\Omega_0,\Ab)\underset{h\to0}{\sim} \frac{|\Omega_0|}{2\pi h},
\end{equation}
where $|\Omega_0|$ denotes the area of $\Omega_0$. The formula carries over for non-simply connected domains but requires imposing Neumann boundary condition on all the boundary components, unlike the setting of the present paper.

Assuming $\lambda\in(0,1)$ is fixed, it is a straightforward  task to adjust the proof in \cite{F} and obtain
\begin{equation}\label{eq:nb-sc*}
N(\lambda h,\Omega,\Ab)\underset{h\to0_+}{\sim} c_1(\lambda)h^{-1/2}\sum_{k \in K^{\Neu}} |\Gamma_k|,\end{equation}
with $c_1(\lambda)$ a positive constant that can be explicitly expressed in terms of $\lambda$. Furthermore, one can adjust the proof in \cite{FK} to obtain
\begin{equation}\label{eq:def-En}
\sum\bigl(\lambda_j(h,\Omega,\Ab)-h \bigr)_-\underset{h\to0_+}{\sim} c_2h^{1/2}\sum_{k\in K^{\Neu}} |\Gamma_k|,
\end{equation}
with $c_2$ a positive explicit constant. 

\subsection{New results for the strip and the annulus}

In both formulas \eqref{eq:nb-sc*} and \eqref{eq:def-En}, we observe that the boundary components with Dirichlet boundary condition are excluded from the asymptotics. This indicates that a similar effect should be observed in the asymptotics of
\[N(h,\Omega,\Ab).\]
The purpose of this paper is to answer this question in the following cases:
\begin{enumerate}[a)]
\item $\Omega=\sS_L:=\mathbb S^1\times(0,L)$ is the strip of radius $1$ and height $L$, and as usual we identify $\mathbb S^1$ and the interval $[0,2\pi)$.
\item $\Omega=\sA_R:=\{x\in \R^2\colon R<|x|<1\}$ is the annulus with inner radius $R$ and outer radius $1$.
\end{enumerate}
Note that, in both these cases, the boundary of $\Omega$ has two components $\Gamma_0$ and $\Gamma_1$. More precisely,
\begin{enumerate}[a$'$)]
\item If $\Omega=S_L$, we set $\Gamma_0=\mathbb S^1\times\{1\}$ and $\Gamma_1=\mathbb S^1\times\{0\}$.
\item If $\Omega=\sA_R$, we set $\Gamma_0=\{|x|=1\}$ and $\Gamma_1=\{|x|=R\}$.
\end{enumerate}
In both cases, we impose Neumann condition on $\Gamma_0$ and Dirichlet condition on $\Gamma_1$, but our analysis would work equally well if the boundary conditions were inversed.

\begin{theorem}\label{thm:strip}
Suppose that $\Omega=\sS_L$ and that the vector potential is
\[\Ab(x_1,x_2)=(-x_2,0).\]
Then, the number of eigenvalues of $\sH_\Omega$ below $h$ satisfies
\[N(h,\Omega,\Ab)\underset{h\to0}{\sim} \frac{|\Omega|}{4\pi h}.  \]
\end{theorem}

We will prove Theorem~\ref{thm:strip} in Section~\ref{sec:strip}. Note that the asymptotics differs drastically from the one in \eqref{eq:nb-sc}. For the strip, imposing Dirichlet boundary condition on one component of the boundary shifts almost half of  the eigenvalues above $h$.
One can intuitively understand this factor of $1/2$ as follows: The projection on the lowest Landau band (without boundary conditions) in ${\mathbb R}^2$ has density of states $(2\pi h)^{-1}$ and consists of exponentially (gaussian!) localized functions. The effect of the Neumann boundary condition will be to slightly lower the energy, whereas the Dirichlet condition will raise it. By the exponential localization, the strongest effect comes from the boundary point closest to a given point in $\Omega$. Thus, the states localized to the half of $\Omega$ closest to $\Gamma_0$ will be pulled down below energy $h$ (and the other half pushed up above).

The scenario  is not exactly the same for the annulus. Imposing Dirichlet boundary condition on the inner boundary of the annulus shifts some of the eigenvalues above $h$, but with a non-trivial fraction.
\begin{theorem}\label{thm:annulus}
Suppose that $\Omega=\sA_R$ and the vector potential is
\[\Ab(x_1,x_2)=\frac12(-x_2,x_1).\]
Then, the number of eigenvalues of $\sH_\Omega$ below $h$ satisfies
\[N(h,\Omega,\Ab)\underset{h\to0}{\sim} \frac1{2h}-\frac{|\sA_R|}{4\pi h|\ln R|} = \frac{1}{2\pi h} |A_{\widetilde{R}}|,  \]
where $\widetilde{R} := \sqrt{\frac{1-R^2}{2|\ln(R)|}} \in (R,\frac{1-R}{2})$.
\end{theorem}
We will prove Theorem~\ref{thm:annulus} in Section~\ref{sec:annulus}.  If we impose Neumann boundary conditions on the inner and outer boundaries of $\sA_R$, then the number of eigenvalues below $h$ is governed by the analogous formula to \eqref{eq:nb-sc}, 
\[N^{\Neu}(h,\sA_R,\Ab)\sim \frac{|\sA_R|}{2\pi h}.\]
This is still consistent with Theorem~\ref{thm:annulus}, since
\[ \forall\,R\in(0,1),\quad R^2<\frac{1-R^2}{2|\ln R|}<1,\]
and there are more eigenvalues below $h$ under the Neumann boundary condition than under mixed Dirichlet-Neumann conditions.

Notice that the simple consideration explained after Theorem~\ref{thm:strip} would lead us to expect that 
$$N(h,\Omega,\Ab) \underset{h\to0}{\sim} \frac{1}{2\pi h} |A_{\frac{1+R}{2}}|,$$
which is different from the result of Theorem~\ref{thm:annulus}. This tells us that the simple intuition is not correct and that the true effect is more subtle. Therefore, understanding the general case would be interesting but we leave that open for future work.

It is worth mentioning the following two remarks on the leading term of $N(h,\Omega,\Ab)$ in Theorem~\ref{thm:annulus}. If $R\to1$, the leading term vanishes and the annulus   $\sA_R$ approaches the unit circle. While if $R\to0$, the leading term is equal to $1/2$ as for the Neumann case in the unit disk (see \eqref{eq:nb-sc}), and $\sA_R$ approaches the unit disk punctured at the origin. Loosely speaking, the last observation suggests that puncturing the unit disk does not have, to leading order, an effect on counting the eigenvalues below the lowest Landau level.

\section{The operator in the strip}\label{sec:strip}

\subsection{Definition of the operator}

Let $L>0$ and consider the infinite strip $\sS_L=\mathbb S^1\times\times(0,L)$. Identifying $\mathbb S^1$ with the interval $[0,2\pi)$, we observe that functions defined on $\sS_L$ must be periodic in the first variable, i.e.
\[u(0,t)=u(2\pi,t)\quad(0<t<L).\]
The magnetic Laplacian we study on $S_L$ is
\[\sH=(-\ii h\partial_s+t)^2-h^2\partial_t^2,\quad  (s,t)\in \sS_L:=\R\times(0,L),\]
with Dirichlet boundary condition on \(t=0\),  Neumann boundary condition on \(t=L\), and periodic conditions with respect to $s$. Decomposing $L^2(\sS_L)$ in Fourier modes with respect to the first variable, the operator $\sH$ is unitary equivalent to the    direct sum of  the fiber operators
\[\fh_m=-h^2\partial_t^2+(t+mh)^2\quad \mbox{on } (0,L),\]
with domain
\[\Dom(\fh_{m})=\{u\in H^2(0,L)\colon  u(0)=0\mbox{ and }u'(L)=0\}.\]
The space \(H^2(0,L)=W^{2,2}(0,L)\) is the Sobolev space of square integrable functions $u$ on $(0,L)$ such that the distributions $u',u''$ are also square integrable on $(0,L)$. The operator $\fh_m$ has compact resolvent and its spectrum is purely discrete consisting of an increasing sequence of simple eigenvalues
\[\lambda_0(\fh_m)<\lambda_1(\fh_m)<\cdots.\]
The spectrum of the operator $\sH$ consists of the \emph{band functions},
\[\spec(\sH)=\overline{\bigcup_{m\in\Z}\spec(\fh^{m})}=\overline{\bigcup_{m\in\Z}\{\lambda_0(\fh_m),\lambda_1(\fh_m),\cdots\}}.\]
As a consequence of the min-max principle and \cite[Proposition~3.2.2]{FH-b}, we have that
\begin{equation}\label{eq:2ndev}
\forall m\in\Z,\quad \lambda_1(\fh_m)>h.
\end{equation}
Consequently, we can identify (by omitting $j=0$ from the notation)
\[\sJ_h:=\{(j,m)\colon \lambda_j(\fh_m)<h\} \approx\{m\in \Z\colon:\lambda_0(\fh_m)<h\}. \]
Denoting by $|\sJ_h|$ the cardinal of of $\sJ_h$, we will prove that as $h\to0$, 
\begin{equation}\label{eq:sL-DN}
|\sJ_h|\sim\frac{L}{2h}.
\end{equation}
The rest of this section is devoted to the proof of \eqref{eq:sL-DN}. Our proof works also for the case where we impose Neumann boundary condition on both boundary components of the strip, for which we find that the  number $\sJ_h^\Neu$ of eigenvalues below $h$ satisfies
\begin{equation}\label{eq:sL-NN}
|\sJ_h^\Neu|\sim \frac{L}{h}.
\end{equation}
\subsection{Harmonic oscillators on the real line}
For $h>0$ and $\xi\in\R$,  consider the operator 
\[T_{\xi,h}=-h^2\partial_t^2+(t+\xi)^2\quad\mbox{on }L^2(\R).\]
The spectrum of $T_h$ consists of the Landau levels $\bigl((2n+1)h\bigr)_{n\geq 0}$ and  an orthonormal basis of corresponding eigenfunctions $(f_n)_{n\geq 0}$ is given as follows
\begin{equation}\label{eq:fn}
 f_n(t)=c_{n,h} H_n\left(\frac{t+\xi}{\sqrt{h}}\right)\exp\left(-\frac{(t+\xi)^2}{2h}\right)\quad(t\in\mathbb R).
\end{equation}
In \eqref{eq:fn}, $c_{n,h}$ is a constant ensuring that $\|f_n\|^2=1$, and $H_n(z)$ is the $n$'th Hermite polynomial. More precisely,
\[c_{n,h}=(2^nn!\sqrt{\pi h})^{-1/2},\quad H_n(z)=(-1)^ne^{z^2}\frac{\dd^n}{\dd z^n}e^{-z^2},\qquad (n=0,1,2\cdots),\]
and we observe, in particular,  that $H_0=1$.  

\subsection{Harmonic oscillators on the half axis}
\subsubsection{Neumann boundary condition.}

For $h>0$ and $\xi\in\R$,  consider the operator 
\[T_{\xi,h}^{\Neu}=-h^2\partial_t^2+(t+\xi)^2\quad\mbox{on }L^2(\R_+),\]
with Neuman boundary condition, $u'(0)=0$, at $t=0$. The eigenvalues of $T^{\Neu}_{\xi,h}$ are all simple and are denoted by,
\[\mu_0^{\Neu}(\xi,h)<\mu_1^{\Neu}(\xi,h)<\cdots.\]
\begin{proposition}\label{prop:N-asy}
As \(\xi/\sqrt{h}\to-\infty\), we have 
\[\mu_0^{\Neu}(\xi,h)-h\sim \pi^{-1/2} h\frac{\xi}{\sqrt{h}}\ee^{-\xi^2/h}.\]
\end{proposition}
Notice that $\xi<0$, so in particular, $h > \mu_0^{\Neu}(\xi,h)$.
\begin{proof}
Consider an orthonormal basis \((u_n)_{n\geq 0}\) of eigenfunctions  such that $u_0>0$. Let us recall the following well known result, which follows from a scaling argument and harmonic approximation (see Prop. 3.2.2 in \cite{FH-b}):
\[\mu_0^{\Neu}(\xi,h)\sim h\qquad(\xi/\sqrt{h}\to -\infty).\]
Harmonic approximation also yields
\begin{equation}\label{eq:app-gs-N}
\| u_0-f_0(\cdot+\xi)\|_{L^2(\R_+)}\to 0\qquad(\xi/\sqrt{h}\to-\infty),
\end{equation}
where $f_0$ is the ground state of the Harmonic oscillator on $\R$.

Moreover, we have from \cite{BH},
\begin{equation}\label{eq:as-gs-N}
u_0(0)\sim \frac{1}{(\pi h)^{1/4}}\ee^{-\xi^2/2h}\qquad(\xi/\sqrt{h}\to-\infty).
\end{equation}

Since \(-h^2f_0''+(t+\xi)^2f_0=h f_0\) and \(-h^2u_0''+(t+\xi)^2u_0=\mu_0^{\Neu}(\xi,h)u_0\), we observe by integration by parts that
\[
\begin{aligned}
h\langle f_0,u_0\rangle_{L^2(\R_+)} &=\langle -h^2f_0''+(t+\xi)^2  f_0,u_0\rangle_{L^2(\R_+)}\\
&=\langle f_0,-h^2u_0''+(t+\xi)^2u_0\rangle_{L^2(\R_+)}+h^2f_0'(0)u_0(0)\\
&=\mu^{\Neu}_0(\xi,h)\langle f_0,u_0\rangle_{L^2(\R_+)}+h^2f_0'(0)u_0(0).
\end{aligned}\]
Consequently, \eqref{eq:as-gs-N} yields,
\[\bigl(h-\mu^{\Neu}_0(\xi,h)\bigr)\langle f_0,u_0\rangle_{L^2(\R_+)}=h^2f_0'(0)u_0(0)\sim -\frac{\xi\sqrt{h}}{\pi^{1/2}}\ee^{-\xi^2/h}.\]
To finish the proof, we use \eqref{eq:app-gs-N}.
\end{proof}

\subsubsection{Dirichlet boundary condition.}

For $h>0$ and $\xi\in\R$,  consider the operator 
\[T_{\xi,h}^{\Dir}=-h^2\partial_t^2+(t+\xi)^2\quad\mbox{on }L^2(\R_+),\]
with Dirichlet boundary condition, $u(0)=0$, at $t=0$. The eigenvalues of $T^{\Dir}_{\xi,h}$ are all simple and are listed as follows,
\[\mu_0^{\Dir}(\xi,h)<\mu_1^{\Dir}(\xi,h)<\cdots.\]
\begin{proposition}\label{prop:D-asy}
As \(\xi/\sqrt{h}\to-\infty\), we have 
\[\mu^{\Dir}_0(\xi,h)-h\sim -\pi^{-1/2} h\frac{\xi}{\sqrt{h}}\ee^{-\xi^2/h}.\]
\end{proposition}
\begin{proof}
By  a scaling argument and harmonic approximation, we have
\[\mu_0^{\Dir}(\xi,h)\sim h\qquad(\xi/\sqrt{h}\to -\infty).\]
Moreover, we have the following rough control of the rate of convergence,
\begin{equation}\label{eq:mu0-Dir}
\mu_0^{\Dir}(\xi,h)- h=\mathcal O\Bigl(\frac1{|\xi|/\sqrt{h}}\Bigr)\qquad(\xi/\sqrt{h}\to -\infty).\end{equation}
By \cite[App. A]{HPRS}, a ground state $u_0$ of \(\mu_0^{\Dir}(\xi,h)\) is necessarily a multiple of the Weber function
\[ U\Bigl(a, \sqrt{2}(t+\xi)/\sqrt{h}\Bigr)\]
where \(a=-\mu_0^{\Dir}(\xi,h)/2h\sim-1/2\) and  \(U(a,-x)\) has the following behavior as \(x\to+\infty\) \cite[Eq.~(A.18)]{HPRS},
\begin{multline*}U(a,-x ) = \frac{\sqrt{2\pi}}{\Gamma(a+\frac12)}\exp(x^2/4)x^{a-\frac12}\Bigl(1+\mathcal O(1/x^2)\Bigr)\\
-\sin(\pi a)\exp(-x^2/4)x^{-a-\frac12}\Bigl(1+\mathcal O(1/x^2)\Bigr).
\end{multline*}
The Dirichlet condition \(u_0(0)=0\) amounts to
\(U(a, \sqrt{2}\xi/\sqrt{h})=0\).  Knowing that the gamma function satisfies \(\Gamma(z)\sim1/z\) as \(z\to0\),  then in the limit \(\xi/\sqrt{h}\to-\infty\), \(a\) is forced to obey the following asymptotics, 
\[ \sqrt{2\pi}\Bigl(a+\frac12\Bigr)\sim \bigl(-\sqrt{2}\xi/\sqrt{h}\bigr)^{-2a}\exp(-\xi^2/h).\]
Using \eqref{eq:mu0-Dir}, we eventually obtain that 
\[a+\frac12\sim -\frac{\xi}{\sqrt{\pi h}}e^{-\xi^2/h}.\]
\end{proof} 

\subsection{Proof of \eqref{eq:sL-DN}}
The idea behind the proof of \eqref{eq:sL-DN} is to localize the angular momenta. Loosely speaking, we will prove that 
\[\lambda_0(\fh_m)<h \]
when the angular momentum satisfies 
\[-L<mh<-L/2,\]
and that $\lambda_0(\fh_m)>h$ otherwise.
Most of the time, it will be more convenient to work with
\begin{equation}\label{eq:def-xi}
\xi=\xi_m:=mh.
\end{equation} 
\subsubsection{Rough localization of the angular momentum}

Let us observe that,   for any $c>1$, if $\xi\geq c\sqrt{h}$ or $\xi<-L-c\sqrt{h}$, then we have
\[\forall t\in[0,L],\quad (t+\xi)^2\geq c^2h>h,\]
hence
\[\lambda_0(\fh_m)>h\mbox{ on }\{m\geq ch^{-1/2}\}\cup\{m\leq -Lh^{-1}-ch^{-1/2}\}.\]
For $\varepsilon\in(0,\frac{L}{4})$, we introduce the intervals,
\[\mathcal I_{\varepsilon}
=\Bigl[-L+\varepsilon,-\frac{L}2-\varepsilon\Bigr],\quad \mathcal J_\varepsilon=\Bigl[-\frac{L}2+\varepsilon,-\varepsilon\Bigr].
\]
To prove \eqref{eq:sL-DN}, it suffices to show that there exists $h_\varepsilon>0$ such that, for $h\in(0,h_\varepsilon]$, we have
\begin{equation}\label{eq:sc-xi}
\lambda_0(\fh_m)<h\mbox{ if }mh\in\mathcal I_\varepsilon,\quad\mbox{and}\quad
\lambda_0(\fh_m)>h\mbox{ if } mh\in\mathcal J_\varepsilon.
\end{equation}
We will prove \eqref{eq:sc-xi} in Propositions~\ref{prop:lb-main} and \ref{prop:lb-main*} below.

\subsection{Favorable angular momenta}

\begin{proposition}\label{prop:lb-main}
Let $\varepsilon\in(0,\frac{L}{4})$. There exists $h_0>0$ such that, if $h\in(0,h_0]$ and $-L+\varepsilon<mh<-\frac{L}2-\varepsilon$, then
\[\lambda_0(\fh_m)<h.\]
\end{proposition}
\begin{proof}
With $\xi=mh$, let $f_0$ be the normalized ground state of the harmonic oscillator, see \eqref{eq:fn}. We introduce the function
\[u=\chi f_0\]
where $\chi$ is a smooth cut-off function satisfying
\[0\leq \chi\leq 1,\quad {\rm supp}\chi\subset[\varepsilon/2,+\infty),\quad \chi|_{[\varepsilon,+\infty)}=1.\]
We compute 
\begin{equation}\label{eq:def-q}
q(u):=\int_{0}^L\bigl( h^2|u'|^2+(t+\xi)^2|u|^2\bigr)\dd t.
\end{equation}
Starting from the identity
\[\bigl\langle T_{\xi,h}f_0,\chi^2f_0\bigr\rangle_{L^2(0,L)}=h\|u\|^2_{L^2(0,L)}\, ,\]
and doing an integration by parts, we get
\[q(u)=h\|u\|_{L^2(0,L)}^2+h^2\|\chi'f_0\|^2_{L^2(0,L)}+h^2f_0'(L)f_0(L).\]
Observe that
\[f_0'(L)f_0(L)=-\frac{L+\xi}{\pi h}\ee^{-(L+\xi)^2/h},\]
and
\begin{multline*}
\|\chi'f_0\|^2\leq\|\chi'\|_\infty^2\int_{{\rm supp}\chi'}|f_0|^2\dd t\leq  
\frac{C_0}{\varepsilon^2}\int_0^\varepsilon \frac{1}{\pi h}\ee^{-(t+\xi)^2/h}\dd t\\
=\mathcal O(\varepsilon^{-2}h^{-1}\ee^{-(\varepsilon+\xi)^2/h}).\end{multline*}
Notice that the condition on $m$ reads as $-L+\varepsilon<\xi<-\frac{L}2-\varepsilon$, hence $|L+\xi|^2<(\varepsilon+\xi)^2$ and
\[q(u)-\|u\|^2_{L^2(0,L)}<0\]
for $h$ sufficiently small. Applying the min-max principle finishes the proof.
\end{proof}

\subsection{Unfavorable angular momenta}

We need the following decay result on the actual ground state of the harmonic oscillator on the half-axis, $T_{\xi,h}^{\Dir}$, with Dirichlet boundary condition.
\begin{proposition}\label{prop:dec}
Let $\alpha\in(0,1)$ and $\varepsilon>0$. There exist constants $C_0,h_0$ such that, if $h\in(0,h_0)$, $mh \in J_\varepsilon$ and
 $u_{m,h}$ is a ground state of $\fh_m$, then 
\[ \int_0^{L} |u_{m,h}|^2\ee^{\alpha(t+mh)^2/h}\dd t\leq C_0\int_0^{L}|u_{m,h}|^2\dd t.\]
\end{proposition}
\begin{proof}
With $\xi=mh<-\varepsilon$, let $\Phi(t)=\frac{\alpha(t+\xi)^2}{2h}$ and $u=u_{m,h}$ be normalized in $L^2$. Given the identity
\[\bigl\langle \fh_m u,\ee^{2\Phi}u\bigr\rangle_{L^2(0,L)}=\lambda_0(\fh_m)\|\ee^{\Phi}u\|^2_{L^2(0,L)},\]
we do an integration by parts. With $q(\cdot)$ from \eqref{eq:def-q}, we get
\[q(\ee^\Phi u)=\lambda_0(\fh_m)\|\ee^\Phi u\|_{L^2(0,L)}^2+h^2\|(\ee^\Phi)' u\|^2_{L^2(0,L)}.\]
Choose $C>0$ such that $(1-\alpha^2)C^2>1$. Decomposing the integrals as follows
\[\int_0^L \cdots\, \dd t=\int_{\{|t+\xi|\leq C\sqrt{h}\}}\cdots \,\dd t +\int_{\{|t+\xi|>C\sqrt{h}\}}\cdots \,\dd t,\]
and using that $|\ee^\Phi|^2$ is uniformly bounded on $\{|t+\xi|\leq C\sqrt{h}\}$ together with the normalization of $u$,  we get
\begin{align*}
    \int_{\{|t+\xi|> C\sqrt{h}\}} h^2|(\ee^\Phi u)'|^2&+(1-\alpha^2)|t+\xi|^2|\ee^\Phi u|^2-\lambda_0(\fh_m)|\ee^\Phi u|^2\,\dd t\\
    &\leq \tilde C \bigl(h+\lambda_0(\fh_m)\bigr),
\end{align*}
where $\tilde C$ is a constant independent of $h$. 

By harmonic approximation, we have $\lambda_0(\fh_m)\sim h$ as $h\to0$, uniformly with respect to $m$ when $mh \in  J_\varepsilon$. Using this  and that the potential term $(1-\alpha^2)|t+\xi|^2$
in the integral is bounded from below by $C^2(1-\alpha^2) h>h$, we get for  $h$ in a sufficiently small right neighborhood of $0$, that
\[\int_{\{|t+\xi|> C\sqrt{h}\}} \Big(|(\ee^\Phi u)'|^2+\frac{h}2\big(C^2(1-\alpha^2)-1\big)|\ee^\Phi u|^2\Big)\dd t\leq \tilde C h.\]
This easily implies the result.
\end{proof}
\begin{proposition}\label{prop:lb-main*}
Let $\varepsilon\in(0,\frac{L}{4})$. There exists $h_0>0$ such that, if $h\in(0,h_0]$ and $mh\in J_\varepsilon$, then
\[\lambda_0(\fh_m)>h.\]
\end{proposition}
\begin{proof}
Let $u$ be the positive normalized ground state of $\fh_m$. We introduce the function
\[f=\chi u\]
where $\chi$ is a smooth cut-off function satisfying
\[0\leq \chi\leq 1,\quad {\rm supp}\chi\subset(-\infty,L-\varepsilon/2),\quad \chi|_{(-\infty,L-\varepsilon]}=1.\]
With $\xi=mh$, we compute  
\begin{equation}\label{eq:def-Q}
Q(f):=\int_{\R_+}h^2 |f'|^2+(t+\xi)^2|f|^2 \,\dd t.
\end{equation}
With $q$ from \eqref{eq:def-q}, we have by the support properties of $f$ and integration by parts,
\[Q(f)=q(f)=\lambda_0(\fh_m)\|\chi u\|_{L^2(0,L)}^2+h^2\|\chi' u\|^2_{L^2(0,L)}.\]
By Proposition~\ref{prop:dec}, for any $\alpha\in(0,1)$, we have
\[
\begin{aligned}
\|f\|^2_{L^2(\R_+)}=\|\chi u\|_{L^2(0,L)}^2&=\|u\|_{L^2(0,L)}^2+\int_{L-\varepsilon}^L(\chi^2-1)|u|^2\dd t\\
&=1+\mathcal O(\ee^{-\alpha|L-\varepsilon+\xi|^2/h}),
\end{aligned}
\]
and
\[\|\chi' u\|^2_{L^2(0,L)}\leq \|\chi'\|_\infty^2\int_{{\rm supp}\chi'}|u|^2\dd t\leq  
\frac{C_0}{\varepsilon^2}\int_{L-\varepsilon}^L|u|^2\dd t\leq \frac{\widetilde{C_0}}{\varepsilon^2}\ee^{-\alpha|L-\varepsilon+\xi|^2/h}).\]
Moreover, by the min-max principle,
\[ Q(f)\geq \mu_0^{\Dir}(\xi,h)\|f\|^2_{L^2(\R_+)}. \]
Collecting the foregoing estimates, we get
\[\lambda_0(\fh_m)\geq \mu_0^{\Dir}(\xi,h)+\mathcal O(\ee^{-\alpha|L-\varepsilon+\xi|^2/h}). \]
To finish the proof, we apply Proposition~\ref{prop:D-asy} and choose $\alpha\in(0,1)$ so that $\alpha|L-\varepsilon+\xi|^2>\xi^2$. In fact, thanks to the condition $-\frac{L}2-\varepsilon<mh=\xi<-\varepsilon$, we have
$\alpha|L-\varepsilon+\xi|^2-\xi^2>\alpha\varepsilon(L-\varepsilon)-(1-\alpha)L^2/4$. It suffices to choose $\alpha\in(0,1)$ such that  $\alpha\varepsilon(L-\varepsilon)-(1-\alpha)L^2/4>0$.
\end{proof}
\section{The operator in the annulus}\label{sec:annulus}
In this section, we fix $R\in(0,1)$ and consider the annulus 
\[\Omega=\sA_R:=\{x\in\R^2\colon R<|x|<1\}.\]
We choose the vector potential as
\[\Ab(x)=\frac12(-x_2,x_1)\]
and consider the operator
\[\sH=(-\ii h\nabla-\Ab)^2\]
in $L^2(\sA_R)$, with   Dirichlet boundary condition on the inner boundary $|x|=R$, and Neumann boundary condition on the outer boundary $|x|=1$.

\subsection{Separation of variables}
 In polar coordinates, the operator $\sH$ can be expressed as 
\[\sH=-\frac{h^2}r\partial_r r\partial_r+\Bigl(\frac{\ii h}{r}\partial_\theta+\frac{r}{2}\Bigr)^2,\]
and if we separate variables, we get the family of operators in $L^2((R,1);r\dd r)$, parameterized by the angular momenta  $m\in\Z$,
\begin{equation}\label{eq:def-annulus-Hm}
\sH_m=h^2\Bigl(-\frac{\dd^2}{\dd r^2}-\frac{1}{r}\frac{\dd}{\dd r}\Bigr)+\Bigl( \frac{mh}{r}-\frac{r}{2}\Bigr)^2, \end{equation}
with Neumann boundary condition at $r=1$ and Dirichlet at $r=R$.

The first eigenvalue of $\sH_m$ is simple with a positive corresponding eigenfunction. We denote the sequence of eigenvalues of $\sH_m$ by
\[\lambda_0(\sH_m)<\lambda_1(\sH_m)\leq \lambda_2(\sH_m)\leq\cdots.\]
With 
\[\sJ_h=\{(j,m)\colon \lambda_j(\sH_m)<h\},\]
and $|\sJ_h|$ denoting the cardinal of $\sJ_h$, we will prove that
\begin{equation}\label{eq:annulus-J}
|\sJ_h|\underset{h\to0}{\sim}\frac1{2h}-\frac{1-R^2}{4 h|\ln R|},
\end{equation}
and this will finish the proof of Theorem~\ref{thm:annulus}.

In the definition of $\sJ_h$, the higher eigenvalues can be excluded. Denoting by $\sH_m^{\Neu}$ the Neumann realization of \eqref{eq:def-annulus-Hm} in $L^2\bigl((0,1),r\dd r\bigr)$, we observe that $\lambda_1(\sH_m)\geq \lambda_1(\sH_m^{\Neu})$. Moreover, $\lambda_1(\sH_m^{\Neu})>h$, as one can see in \cite[Prop.~2.1]{FFGKS}. Consequently, $\sJ_h$ reduces to
\begin{equation}\label{eq:annulus-J-def*}
\sJ_h=\{m\in\Z\colon \lambda_0(\sH_m)<h\}.
\end{equation}

\begin{remark}
The result of Theorem~\ref{thm:annulus} can be understood through a simple picture in radial coordinates. The radial eigenfunction (non-normalized) to the magnetic operator with constant magnetic field on all of ${\mathbb R}^2$ and angular momentum $m$ is given by $u(r) = \exp(-\phi(r)/4h)$ with $\phi(r) = r^2 - 2 r_{*}^2 \ln(r)$, with $r_{*} = \sqrt{2mh}$ (see \eqref{eq:qm-disc*}, \eqref{eq:def-r*}, \eqref{eq:def-phi} below) and is strongly localized near $r= r_{*}$.
The (positive) effect of the Dirichlet boundary condition at $r=R$ will be governed by the decay $\exp(-\phi(r=R)/4h)$ and the (negative) effect of the Neumann boundary condition at $r=1$ by $\exp(-\phi(r=1)/4h)$. Thus, the effect of the Neumann boundary condition will be strongest, and pull the energy down, exactly for those $m$ where $\phi(1) < \phi(R)$. This line of thought leads to the introduction of the sets of momenta corresponding to $\mathcal I_{\varepsilon,h}$ and $\mathcal J_{\varepsilon,h}$ below and to the strategy of the proof.
\end{remark}

\subsection{Rough localization of angular momentum}
For any $c>1$, if
\[m\geq \frac1{2h}+\frac{c}{\sqrt{h}},\]
then, for all $r \in (R,1)$,
\[v_m(r):=\Bigl( \frac{mh}{r}-\frac{r}{2}\Bigr)^2\geq h^2\Bigl(m-\frac{r^2}{2h}\Bigr)^2\geq c^2h>h,\]
and $\lambda_0(\sH_m)\geq c^2h>h$.

Similarly, if $c > R$ and 
\[m\leq \frac{R^2}{2h}-\frac{c}{\sqrt{h}},\]
then
\[v_m(r)\geq\frac{c^2h}{R^2}>h,\]
and $\lambda_0(\sH_m)\geq c^2h>h$.

Therefore, for every $\varepsilon > 0$, there is $h_0>0$ such that 
\[\sJ_h\subset \Bigl(\frac{R^2}{2h}-\frac{R+\varepsilon}{\sqrt{h}}, \frac1{2h}+\frac{1+\varepsilon}{\sqrt{h}}\Bigr),\quad\forall\,h\in(0,h_0].\]

For every $\varepsilon<\frac12\min\bigl(1-\frac{1-R^2}{2|\ln R|},\frac{1-R^2}{2|\ln R|}-R^2\bigr)$ and $h>0$, we introduce the intervals 
\[\mathcal I_{\varepsilon,h}=\Biggl[\frac{1-R^2}{4h|\ln R|}+\frac{\varepsilon}{2h}, \frac{1-\varepsilon}{2h}\Biggr]\quad\mbox{ and }\quad\mathcal J_{\varepsilon,h}=\Biggl[\frac{R^2+\varepsilon}{2h}, \frac{1-R^2}{4h|\ln R|}-\frac{\varepsilon}{2h}\Biggr].\]
To prove \eqref{eq:annulus-J}, it suffices to show that there is $h_\varepsilon>0$ such that, for $h\in(0,h_\varepsilon]$, we have
\begin{equation}\label{eq:annulus-J*}
\mathcal I_{\varepsilon,h} \cap {\mathbb Z}\subset \sJ_h\quad\mbox{ and }\quad
\mathcal J_{\varepsilon,h}\cap \sJ_h=\emptyset.
\end{equation}
We will prove \eqref{eq:annulus-J*} in Propositions~\ref{prop:lb-disc} and \ref{prop:lb-disc*} below. The proof essentially follows the same  scheme as the one in Section~\ref{sec:strip}, with additional difficulties related to the non-flat geometry in this model.%
\subsection{A family of quasi modes}
For $m\in\Z$, the function
\begin{equation}\label{eq:qm-disc}
u(r)=r^m\exp(-r^2/4h)
\end{equation}
satisfies $\sH_m u=hu$ on $\R_+$. Furthermore, when $m\in \mathcal I_{\varepsilon,h}\cup \mathcal J_{\varepsilon,h}$,
$u$ is localized at 
\begin{equation}\label{eq:def-r*}
r_* = r_{*}(m):=\sqrt{2mh}\in\Bigl[\sqrt{R^2+\varepsilon},\sqrt{1-\varepsilon}\,\Bigr].
\end{equation}
 To see this, we write 
\begin{equation}\label{eq:qm-disc*}
u(r)=\exp\bigl(-\phi(r)/4h\bigr),
\end{equation}
where 
\begin{equation}\label{eq:def-phi}
\phi(r)=r^2-2r_*^2\ln r.
\end{equation} 
Then, assuming \eqref{eq:def-r*},  we have
\[\min_{r\in [R,1]}\phi(r)=\phi(r_*)=r_*(1-2\ln r_*)>0,\quad \phi''(r_*)=4,\]
and by Laplace's method,
\[\int_{R}^1|u(r)|^2r\dd r\underset{h\to0}{\sim} r_*\sqrt{\frac{2\pi}{\phi''(r_*)}}\ee^{-2\phi(r_*)/4h}=\sqrt{\pi/4}\,r_*\ee^{-\phi(r_*)/2h}.\]
Loosely speaking, we will show that the angular momenta $m$ corresponding to $r_*$ such that $R^2<r_*^2<\frac{1-R^2}{2|\ln R|}$  yield ground state energies below $h$, while for angular momenta corresponding to $\frac{1-R^2}{2|\ln R|}<r_*^2<1$, we will show that the ground state is above $h$. That is analogous to the  observation  in the flat case (Section~\ref{sec:strip}).
\subsection{Favorable angular momenta}

\begin{proposition}\label{prop:lb-disc}
Given $0<R<1$ and $0<\varepsilon<\frac12\bigl(1-\frac{1-R^2}{2|\ln R|}\bigr)$,  there exists $h_0>0$ such that, if $h\in(0,h_0]$ and $\frac{1-R^2}{2|\ln R|}+\varepsilon<2mh<1-\varepsilon$, then
\[\lambda_0(\sH_m)<h.\]
\end{proposition}
\begin{proof}
Note that $r_*=\sqrt{2mh}$ satisfies
\[R^2+\varepsilon<\frac{1-R^2}{2|\ln R|}+\varepsilon<r_*^2<1-\varepsilon\]
and that $\phi(1)<\phi(R)-2\varepsilon|\ln R|<\phi(R)$, where $\phi$ is introduced in \eqref{eq:def-phi}. Choose a positive constant $\eta<\varepsilon$ such that $R+\eta<r_*$ and
\[\phi(1)<\phi(R+\eta)-\varepsilon|\ln R|<\phi(R+\eta).\]
Moreover, since $\phi'\leq 0$ on $(0,r_*]$,  we have $\phi(R+\eta)\leq \phi(r)\leq \phi(R)$ for all $r\in[R,R+\eta]$.

We introduce the function
\[f=\chi u,\]
where $u$ is the function in \eqref{eq:qm-disc}, and  $\chi$ is a cut-off function satisfying
\[0\leq \chi\leq 1,\quad {\rm supp}\chi\subset(R,+\infty),\quad \chi|_{[R+\eta,+\infty)}=1.\]
Starting from the identity
\[\int_R^1\sH_m u \,\chi^2u\,\dd r=h\int_R^1 |\chi u|^2r\dd r,\]
and doing an integration by parts, we get
\begin{align}
q_m(f)&:=\int_R^1\left( h^2|f'|^2+\Bigl( \frac{mh}{r}-\frac{r}{2}\Bigr)^2|f|^2\right)r\dd r\label{eq:def-qm-disc}\\
&=h\int_R^1 |f|^2r\dd r+\int_R^1|\chi'u|^2r\dd r+u'(1)u(1).\nonumber
\end{align}
Thanks to \eqref{eq:qm-disc*} and $R+ \eta < r_{*}$, we have on the one hand,
\[\int_R^1|\chi'u|^2r\dd r\leq \frac{C}{\eta^2}\int_R^{R+\eta}|u|^2r\dd r\leq \frac{C}{\eta}\ee^{-\phi(R+\eta)/2h}.\]
On the other hand, we have,
\[u'(1)u(1)=\Bigl(m-\frac1{2h}\Bigr)\ee^{-1/2h}= \frac{r_*^2-1}{2h}\ee^{-\phi(1)/2h}.\]
Consequently,
\[\int_R^1|\chi'u|^2r\dd r+u'(1)u(1)\leq \frac{C}{\eta}\ee^{-\phi(R+\eta)/2h} +\frac{r_*^2-1}{2h}\ee^{-\phi(1)/2h}.\]
Since $\phi(1)<\phi(R+\eta)-\varepsilon|\ln R|$, we get for $h$ small enough,
\[\frac{C}{\eta}\ee^{-\phi(R+\eta)/2h} +\frac{r_*^2-1}{2h}\ee^{-\phi(1)/2h}<0,\]
and therefore,
\[ q_m(f)-h\int_R^1 |f|^2r\dd r<0 .\]
Applying the min-max principle, we get that $\lambda_0(\sH_m)<h$.
\end{proof}

\subsection{Unfavorable angular momenta}
\begin{proposition}\label{prop:lb-disc*}
Given $\varepsilon\in(0,\frac{1-R^2}{4})$, there exists $h_0>0$ such that, if $h\in(0,h_0]$ and $R^2+\varepsilon<2mh<\frac{1-R^2}{2|\ln R|}-\varepsilon$, then
\[\lambda_0(\sH_m)> h.\]
\end{proposition}
To prove Proposition~\ref{prop:lb-disc*} we need a result on the decay of the true ground states.
\begin{proposition}\label{prop:dec-disc}
Given $0<R<1$, $0<\varepsilon<(1-R)/2$ and $\alpha\in(0,1)$, there exist constants $C_0,h_0$ such that, if $h\in(0,h_0]$, $r_*:=\sqrt{2mh}\in[R+\varepsilon,1-\varepsilon]$ and
 $u_{m}$ is a ground state of $\sH_m$, then 
\[ \int_R^1 |u_{m}|^2\ee^{\alpha\frac{\phi(r)-\phi(r_*)}{2h}}\dd r\leq C_0\int_R^1|u_{m}|^2\dd r,\]
where
\[\phi(r)=r^2-2r_*^2\ln r.\]
\end{proposition}
\begin{proof}
Similar to the proof of Proposition~\ref{prop:dec}. Firstly, notice that the potential
\[v_m(r):=\left(\frac{mh}{r}-\frac{r}2\right)^2\]
has a unique minimum at $r_*$ with $v_m(r_*)=0$ and $v_m''(r_*)=2$, and by harmonic approximation, we have
\[\lambda_0(\sH_m)\underset{h\to0}{\sim} h,\quad \lambda_1(\sH_m)\underset{h\to0}{\sim} 3h.\]
 Secondly, we introduce the function
\[ \Phi(r)=\alpha\frac{\phi(r)-\phi(r_*)}{4h}\]
and notice that it is minimal for $r=r_*$, and that
\[ h^2|\Phi'(r)|^2=\alpha^2v_m(r).\] 
In the sequel, $u$ denotes the positive and normalized ground state of $\sH_m$, and we write $\lambda$ for $\lambda_0(\sH_m)$, the ground state energy for $\sH_m$. Starting from the identity
\[\bigl\langle \sH_mu,\ee^{2\Phi}u\bigr\rangle=\lambda\|\ee^\Phi u\|^2\]
and doing an integration by parts, we get, with $q_m(\cdot)$  introduced in \eqref{eq:def-qm-disc},
\[q_m(\ee^\Phi u)=\lambda\|\ee^\Phi u\|^2+h^2\|(\ee^\Phi)' u\|^2.\]
Choose $C>0$ such that $(1-\alpha^2)C^2>16$. Decomposing the integrals as follows
\[\int_R^1 \cdots\, r\dd r=\int_{\phi'(r)\leq C\sqrt{h}}\cdots \,r\dd r +\int_{\phi'(r)>C\sqrt{h}}\cdots \,r\dd r\]
and using that the integral of $|\ee^\Phi u|^2$ over $\{\phi'(r)\leq C\sqrt{h}\}$ is uniformly bounded (thanks to the mean value theorem),   we get
\[\int_{\phi'(r)> C\sqrt{h}} \Big(|(\ee^\Phi u)'|^2+(1-\alpha^2)v_m(r)|\ee^\Phi u|^2-\lambda|\ee^\Phi u|^2\Big)\dd t\leq \tilde Ch,\]
where $\tilde C$ is a constant independent of $h$. 

Using that $\lambda\sim h$ and that the potential term $(1-\alpha^2)v_m(r)$
in the integral is bounded from below by $ C^2(1-\alpha^2)h/16>h$,
we get for   $h$ in a sufficiently small right neighborhood of $0$, that
\[\int_{\phi'(r)> C\sqrt{h}} \Big(|(\ee^\Phi u)'|^2+{\frac{h}2}\big(\mbox{$\frac{1}{16}$}C^2(1-\alpha^2)-1\big)|\ee^\Phi u|^2\Big)\dd t\leq \tilde Ch.\]
\end{proof}
Another ingredient of the proof of Proposition~\ref{prop:lb-disc*} is the behavior of the Dirichlet ground state energy on a semi-axis (which corresponds to the exterior of the disc with radius $R$). Let $\sH_m^{\Dir}$ denote the operator $\sH_m$ on $(R,+\infty)$ with Dirichlet boundary condition at $r=R$.
\begin{proposition}\label{prop:disc-D-asy}
Given $R,\varepsilon>0$, there exist positive constants  $c_*<c^*$ and $h_0$ such that, if 
\[m\in\N,\quad  h\in(0,h_0], \quad r_*:=\sqrt{2mh}>R+\varepsilon,\] then we have 
\[h+c_*Rh^{1/2}\ee^{-c(r_*)/h}\leq \lambda(\sH_m^{\Dir})\leq h+c^*Rh^{1/2}\ee^{-\frac{\phi(R)-\phi(r_*)}{2h}},
\]
where
\[\phi(r)=r^2-2r_*^2\ln r.  \]
\end{proposition}
\begin{proof}
By scaling, it suffices to prove the proposition for $R=1$ and $1+\varepsilon <\sqrt{2mh}<C$ with $C$ a fixed constant.  We write $\lambda$ and $\lambda_1$ for $\lambda(\sH_m^{\Dir})$ and $\lambda_1(\sH_m^{\Dir})$ respectively, and note that $\lambda\sim h$, and $\lambda_1\sim 3h$,  by harmonic approximation. Moreover, by domain monotonicity and the min-max principle, we may compare with the harmonic oscillator on  $\mathbb R$ and get that $\lambda>h$. 

The function $u(r)=\ee^{-\phi(r)/4h}=r^m\ee^{-r^2/4h}$ introduced  in \eqref{eq:qm-disc} satisfies
\[\sH_mu=hu,\]
and by Laplace's approximation,
\[
 M_h:=\int_{1}^{+\infty} |u(r)|^2\, r\dd r\underset{h\to0}{\sim}\sqrt{\pi h } r_*\ee^{-\phi(r_*)/2h}.
 \]
We introduce the function
\[v=f u,\quad f(r)=\int_1^r\chi(\rho)\rho^{-1-2m}\ee^{\rho^2/2h}\dd\rho, \]
where $\chi\in C_c^\infty(\R;[0,1])$ satisfies $\chi=1$ on $(-\infty,1+\frac12r_*]$ and $\chi=0$ on $[1+r_*,+\infty)$. We have the uniform bound
\[ 0\leq f(r)\leq \int_1^{1+r_*}\rho^{-1-2m}\ee^{\rho^2/2h}\dd\rho=\mathcal O(h\ee^{1/2h}).\]
Moreover, if $r\geq 1+2h$,  we can bound $f(r)$ from below as
\[f(r) \geq \int_1^{1+2h}\rho^{-1-2m}\ee^{\rho^2/2h}\dd\rho\geq c_0h\ee^{1/2h},\]
where $c_0>0$ is independent from $h\in(0,1]$.

Notice that $v$ satisfies the Dirichlet condition, $v(1)=0$, and
\begin{equation}\label{eq:disc-Hv}
(\sH_m-h)v=-h^2\chi'(r)r^{-m-1}\ee^{r^2/4h}.
\end{equation}
We calculate 
\[\begin{aligned}
\langle (\sH_m-h)v,v\rangle&=-h^2\int_1^{+\infty}\chi'(r)f(r)\dd r\\
&=h^2\int_1^{+\infty}\chi(r)^2r^{-1-2m}\ee^{r^2/2h}\dd r\geq c_1h^3\ee^{1/2h},\\
\|(\sH_m-h)v\|^2&=h^4\int_{1}^{+\infty} \Bigl(\chi'(r)\Bigr)^2 r^{-2m-1}\ee^{r^2/2h}\dd r\\
&\leq c_2 h^5\ee^{1/2h}.
\end{aligned} \]
Similarly,  we check that
\[ c_0^2h^2\ee^{1/h}M_h\leq \|v\|^2=\int_1^{+\infty}|f(r)|^2|u(r)|^2r\dd r=\mathcal O(h^2\ee^{1/h}M_h),\]
and
\[ \langle (\sH_m-h)v,v\rangle =\mathcal O\bigl(h^3\ee^{1/2h}\bigr).\]
Now we apply Temple's inequality. With 
\[\eta=\frac{\langle (\sH_m-h)v,v\rangle}{\|v\|^2},\quad 
\epsilon^2=\frac{\|(\sH_m-h)v\|^2}{\|v\|^2}-\eta^2\mbox{ and }\beta=h<\lambda_1-h,\]
 we have
\[ \lambda-h \geq \eta-\frac{\epsilon^2}{\beta-\eta}\geq c_*h^{1/2}\bigl(1+o(1)\bigr)\ee^{\frac{\phi(r_*)-1}{2h}},\]
where $c_*$ is a positive constant, and
\[\phi(r_*)-1=\phi(r_*)-\phi(1)<0\]
because $r_*>1$. Moreover, by the mix-max principle, we have
\[\lambda-h\leq \eta\leq c^*h^{1/2}\bigl(1+o(1)\bigr)\ee^{\frac{\phi(r_*)-1}{2h}}.\]
\end{proof}

\begin{proof}[Proof of Proposition~\ref{prop:lb-disc*}]
The condition on $m$ reads as
\[R^2+\varepsilon<r_*^2<\frac{1-R^2}{2|\ln R|}-\varepsilon<1-\varepsilon.\]
Consequently, we have
\[ \phi(1)>\phi(R)+2\varepsilon|\ln R|.\]
Choose a positive $\eta<\varepsilon$ such that $1-\eta>r_*$ and
\begin{equation}\label{eq:cond-eta-disc}
\phi(1-\eta)>\phi(R)+\varepsilon|\ln R|.
\end{equation}
Let $\psi$ be a normalized  ground state of $\sH_m$. We introduce the function
\[f=\chi \psi\]
where $\chi$ is a cut-off function satisfying
\[0\leq \chi\leq 1,\quad {\rm supp}\chi\subset(-\infty,1-\eta/2),\quad \chi|_{(-\infty,1-\eta]}=1.\]
We compute  
$q_m(f)$ where $q_m(\cdot)$ is introduced in \eqref{eq:def-qm-disc}. After integration by parts, we have
\[q_m(f)=\lambda(\sH_m)\|\chi \phi\|^2+h^2\|\chi' \phi\|^2.\]
By Proposition~\ref{prop:dec-disc}, for any $\alpha\in(0,1)$, we have
\[
\begin{aligned}
\|f\|^2=\|\chi \psi\|^2&=\|\psi\|^2+\int_{1-\eta}^1(\chi^2-1)|\psi|^2r\dd r\\
&=1+\mathcal O\Bigl(\ee^{-\alpha\frac{\phi(1-\eta)-\phi(r_*)}{2h}}\Bigr),
\end{aligned}
\]
and
\[\|\chi' \psi\|^2\leq \|\chi'\|_\infty\int_{{\rm supp}\chi'}|u|^2r\dd r\leq  
\frac{C_0}{\eta^2}\int_{1-\eta}^1|\psi|^2r\dd r=\mathcal O\Bigl(\ee^{-\alpha\frac{\phi(1-\eta)-\phi(r_*)}{2h}}\Bigr).\]
Moreover, since $f$ is supported in $(R,+\infty)$, we have by the min-max principle and Proposition~\ref{prop:disc-D-asy}, 
\[ q_m(f)\geq \lambda(\sH_m^{\Dir})\|f\|^2\geq \Bigl(h+c_*Rh^{1/2}\ee^{-\frac{\phi(R)-\phi(r_*)}{2h}}\Bigr)\|f\|^2. \]
Collecting the previous estimates, we get
\[\lambda(\sH_m)\geq h+c_*Rh^{1/2}\ee^{-\frac{\phi(R)-\phi(r_*)}{2h}}+\mathcal O\Bigl(\ee^{-\alpha\frac{\phi(1-\eta)-\phi(r_*)}{2h}}\Bigr). \]
To finish the proof, we  choose $\alpha\in(0,1)$ so that 
$\alpha\bigl(\phi(1-\eta)-\phi(r_*)\bigr)>\phi(R)-\phi(r_*)$,  which is possible, thanks to \eqref{eq:cond-eta-disc}.
\end{proof}

\subsection*{Acknowledgments}{\small 
SF was partially supported by the Villum Centre of Excellence for the Mathematics of Quantum Theory (QMATH) with Grant No.10059. SF was partially funded by the European Union. Views and opinions expressed are however those of the authors only and do not necessarily reflect those of the European Union or the European Research Council. Neither the European Union nor the granting authority can be held responsible for them.\\
AK is partially supported by CUHK-SZ grant no. UDF01003322.
}


\begin{thebibliography}{10}
%
\bibitem{BH}
C. Bolley,  B. Helffer. An application of semi-classical analysis to the asymptotic study of the supercooling field of
a superconducting material.
Annales de l'I. H. P., section A, tome 58, no 2 (1993), p. 189-233.

\bibitem{FFGKS} S. Fournais, R.L. Frank, M. Goffeng, A. Kachmar, M.P. Sundqvist.
Counting negative eigenvalues for the magnetic Pauli operator. arXiv:2308.14680 (2023), accepted in Duke Math. J.
%
%
\bibitem{FH-b} S. Fournais, B. Helffer. Spectral methods in surface superconductivity.
Progr. Nonlinear Differential Equations Appl., 77
Birkhäuser Boston, Inc., Boston, MA, 2010.
%
\bibitem{FK} S. Fournais, A. Kachmar.
On the energy of bound states for magnetic Schr\"odinger operators.
J. Lond. Math. Soc. (2), 80 (2009), 233--255.
%
\bibitem{F} R. L. Frank. On the asymptotic number of edge states for magnetic Schr\"odinger operators. Proc. London Math. Soc. (3), 95 (2007) 1--19. 
%
\bibitem{HPRS} P.D. Hislop, N. Popoff,  N. Raymond, M.P. Sundqvist.
Band functions in the presence of magnetic steps. 
Math. Models Methods Appl. Sci.26 (2016), no.1, 161-184.
\end{thebibliography}
\end{document}